\documentclass[10pt]{article}

\usepackage{textcomp}
\usepackage{cite}      
\usepackage{graphicx}  
\usepackage{amsmath}
\usepackage{authblk}
\usepackage[tight,TABTOPCAP]{subfigure}
\usepackage[margin=1in]{geometry}
\usepackage[colorlinks=true, urlcolor=blue]{hyperref}

\usepackage[usenames,dvipsnames]{xcolor}
\definecolor{blue}{rgb}{0.235,0.310,1}
\definecolor{gray}{rgb}{0.9,0.9,0.9}
\definecolor{green}{rgb}{0.0,0.635,0.067}
\definecolor{violet}{rgb}{0.761,0.0,0.984}

\usepackage{listings}
\newcommand{\de}{\,\mathrm{d}} 
\newcommand{\units}[1]{~\mathrm{#1}} 
\newcommand{\figref}[1]{Fig.~\ref{#1}} 
\newcommand{\tabref}[1]{Table~\ref{#1}} 

\graphicspath{{img/}}

\title{Errata and comments for ``Numerical and analytical modeling of busbar systems''}

\author{Luca Giaccone\thanks{corresponding author: luca.giaccone@polito.it} }
\author{Aldo Canova}
\affil{\small Politecnico di Torino, Dipartimento Energia}

\date{\today}

\begin{document}
\maketitle

\begin{abstract}
This note covers two parts. The first one provides an errata to the paper ``Numerical and analytical modeling of busbar systems''. We mainly give the correction for three equations affected by a typographical mistake. Despite the corrections that we are providing with this note, we think that the implementation of these equations can be quite onerous. Hence, in the second part of this document we provide the download link to our implementation of the equations (developed in MATLAB environment). Moreover, to help in using these functions, we explain their behavior by means of some examples.
\end{abstract}

\section*{Acknowledgement}
We would like to thank John Compter for contacting us about reference \cite{Canova2009}. His comments allowed the authors to identify the typographical errors that we are going to correct in this note.

\section{Errata Corrige}
\begin{figure}[!ht]
   \centering
		\includegraphics[angle=0, width=0.4\columnwidth]{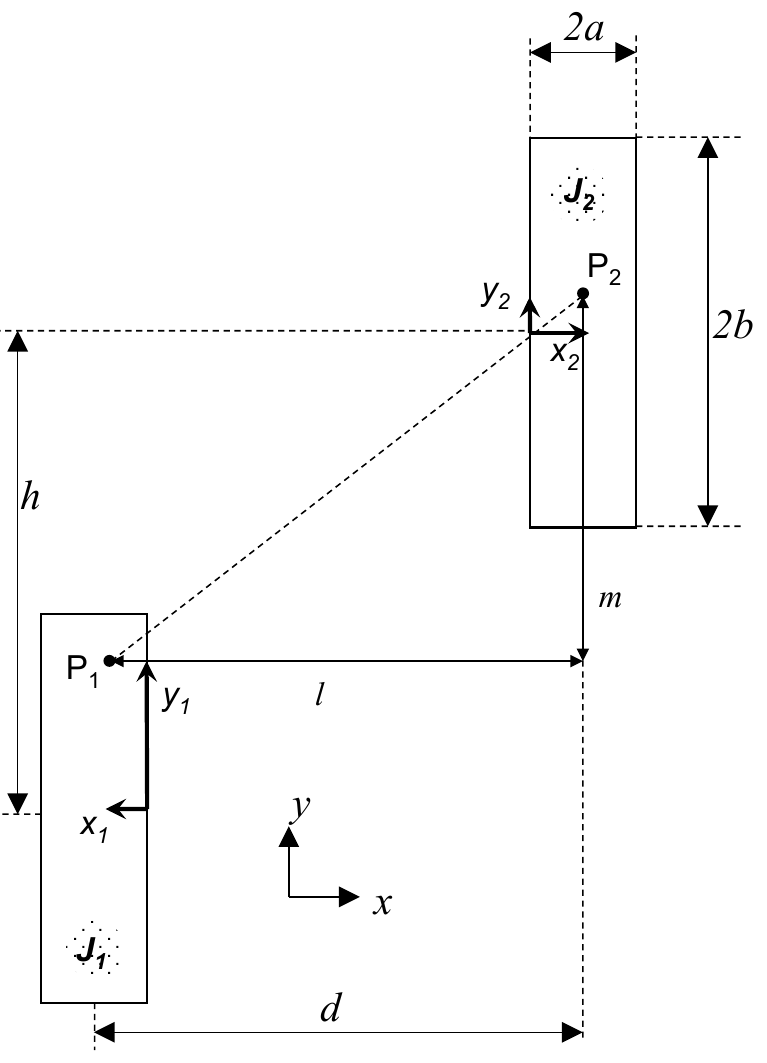}
   \caption{Layout of non-adjacent massive conductors.}
	  \label{fig:non_adjacent}
\end{figure}

The paper ``Numerical and analytical modeling of busbar systems'' deals with the computation of electrodynamic forces between massive conductors. In section III.B we propose an extension for the IEC865-93 \cite{IEC865-93} that makes possible the computation of the forces between non adjacent conductors as the one described in \figref{fig:non_adjacent}. With reference to the appendix of \cite{Canova2009}, the corrections summarized in \tabref{tab:corrections} should be done to have proper results.
\begin{table}[!ht]
\centering
\renewcommand{\arraystretch}{2.5} 
\caption{Table of corrections}
\vspace{0.1cm}
\begin{tabular}{c | c}
\hline
\hline
{\bf Errata} & {\bf Corrige} \\
\hline
\hline
$S=8b^3-12b^2$ & $S=8b^3-12b^2\textcolor{blue}{h}$ \\
\hline
$T=8b^3+12b^2$ & $T=8b^3+12b^2\textcolor{blue}{h}$ \\
\hline
$\begin{array}{l}
F_{x} =  \bigg[ 2A(I+M+G) \arctan\left(\dfrac{C}{A}\right) \\
  + 2B\left(I-N+G\right) \arctan\left(\dfrac{C}{B}\right)\\
  -4h\left(I+G\right) \arctan\left(\dfrac{C}{\textcolor{red}{H}}\right) \ldots\\
  \ldots \ldots \ldots
\end{array}$ &
$\begin{array}{lll}
F_{x}  =  \bigg[ 2A(I+M+G) \arctan\left(\dfrac{C}{A}\right) \\
  + 2B\left(I-N+G\right) \arctan\left(\dfrac{C}{B}\right)\\
  -4h\left(I+G\right) \arctan\left(\dfrac{C}{\textcolor{blue}{h}}\right) \ldots\\
  \ldots \ldots \ldots
\end{array}$ \\
\hline
inset of Fig. 14: & inset of Fig. 14: \\
\includegraphics[scale=0.8]{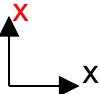} & \includegraphics[scale=0.8]{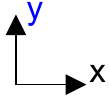} \\
\hline
\hline
\end{tabular}
\label{tab:corrections}
\end{table}%
\subsection{Other remarks}
In addition to the errata corrige, the authors would like to provide some other remarks to avoid further confusion about \cite{Canova2009}.

\noindent
{\bf Comment about equations (22), (23) and (24):} all the equations describe a four fold integral. Some confusion can be caused by the fact that the limits of integration appear in a sequence that does not match the corresponding differentials. Hence, keeping the same numeration of the original paper, these equations should be more conveniently written as:
\begin{equation}
\label{eq:force_adj}
F_{x}=\int^{2a}_{0} \!\!\! \int^{2a}_{0} \!\!\! \int^{b}_{-b} \! \int^{b}_{-b} f \dfrac{l}{\sqrt{l^2+m^2}} \de x_1 \de x_2 \de y_1 \de y_2 \tag{22}
\end{equation}
\begin{equation}
\label{eq:force_nadj_h}
F_{x}=\int^{2a}_{0} \!\!\! \int^{2a}_{0} \!\!\! \int^{b}_{-b} \! \int^{b}_{-b} f \dfrac{l}{\sqrt{l^2+m^2}} \de x_1 \de x_2 \de y_1 \de y_2 \tag{23}
\end{equation}
\begin{equation}
\label{eq:force_nadj_v}
F_{y}=\int^{2a}_{0} \!\!\! \int^{2a}_{0} \!\!\! \int^{b}_{-b} \! \int^{b}_{-b} f \dfrac{m}{\sqrt{l^2+m^2}} \de x_1 \de x_2 \de y_1 \de y_2 \tag{24}
\end{equation}

\noindent
{\bf Comment on the range of validity for the equations related to non-adjacent conductors:} it is worth remarking that, with reference to \figref{fig:non_adjacent}, the proposed equations for the force computation are valid only if the following constraints are satisfied:
\begin{itemize}
\item $d > 2a$ (please note that we use $>$ and not $\ge$)
\item $h > 2b$ (please note that we use $>$ and not $\ge$)
\end{itemize}

\noindent
{\bf Comment about the current density inside massive conductors:} the proposed equations are an extension of what can be found in the standard IEC865-93 \cite{IEC865-93} where the current density inside each massive conductor is assumed to be homogeneous. Only making this assumption it is possible to develop the proposed analytical equations. Hence, for all the cases where skin and proximity effect are not negligible, the proposed equations must be used together with the correction factors that are also provided in \cite{Canova2009}.

\section{MATLAB functions and their use}
Even if we did our best to shorten the length of the formulas proposed in \cite{Canova2009}, we believe that the implementation of these analytical equations can be quite onerous. We are writing this document because we did three mistakes in writing them, therefore, it is possible to make other mistakes trying to implement the equations. For these reasons, we would like to share our codes with anyone is interested. The MATLAB functions can be downloaded by this website \href{http://www.cadema.polito.it/}{{\tt www.cadema.polito.it}} visiting the download section related to the analytical methods (\href{http://www.cadema.polito.it/?page_id=508}{direct link}).

\begin{figure}[!ht]
\centering
\subfigure[]{\includegraphics[scale=0.41]{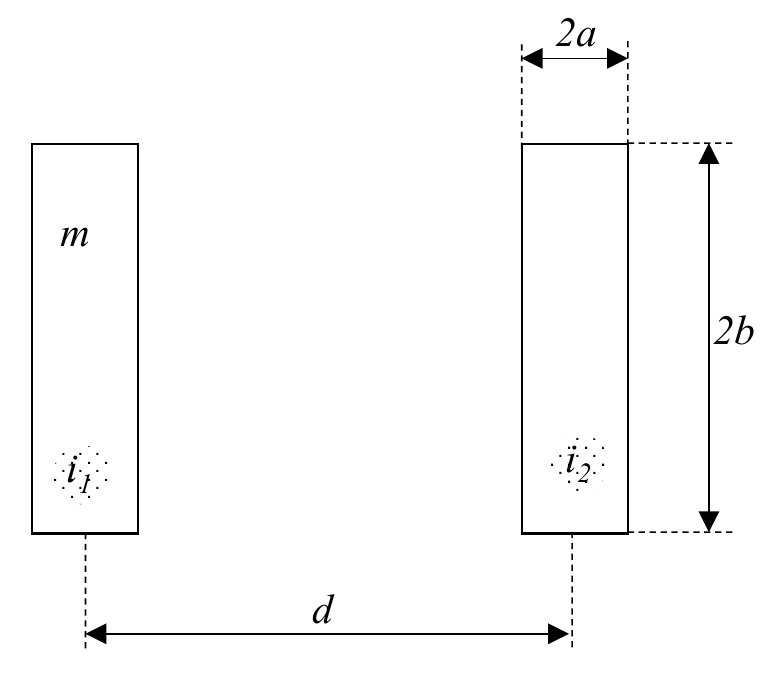}\label{fig:parametrical_adjacent}} \par
\subfigure[]{\includegraphics[scale=0.41]{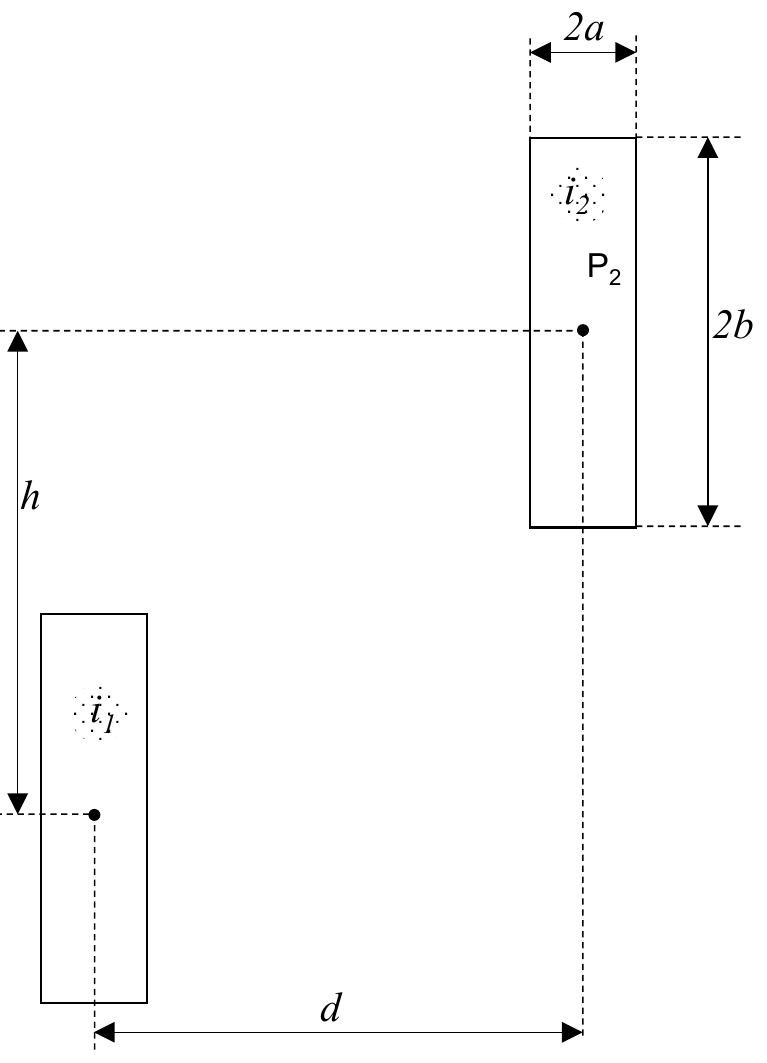}\label{fig:parametrical_non_adjacent}} %
\caption{Reference configuration for adjacent conductors (a). Reference configuration for non-adjacent conductors (b).}
\end{figure}

In the following pages three examples are presented to explain how to use properly the MATLAB functions. The first two examples are related to adjacent conductors as the ones in \figref{fig:parametrical_adjacent}. The third example is related to non-adjacent conductors as the ones in \figref{fig:parametrical_non_adjacent}. In all the examples the behavior of the MATLAB functions is explained describing a code list that can be easily reproduced. Moreover, the results of each example is further validated by means of a 2D Finite Element Method (FEM). It is worth noting that we are using an independent FEM code \cite{Meeker2006} to validate our analytical equations.

%
\subsection{Example \#1: adjacent conductors}
This example describes the use of the function {\tt adjacent.m}. It is carried out a parametrical analysis by varying the parameter $d$.

\noindent
The x-force can be computed by means of the following code. 
\begin{lstlisting}
%% Section 1: Geometrical parameters
a = 0.005;                  % x dimension (m)
b = 0.05;                   % y dimension (m)
d = linspace(0.011,0.2,15); % distance between conductors (m)

%%  Section 2: conductors currents
i1 = 1; % current in conductor 1 (A)
i2 = 1; % current in conductor 2 (A)

%%  Section 3: computation
for i = 1:length(d)
    Fx(i) = adjacent(a,b,d(i),i1,i2);
end
\end{lstlisting}

\begin{figure}[!ht]
   \centering
		\includegraphics[angle=0, width=1\columnwidth]{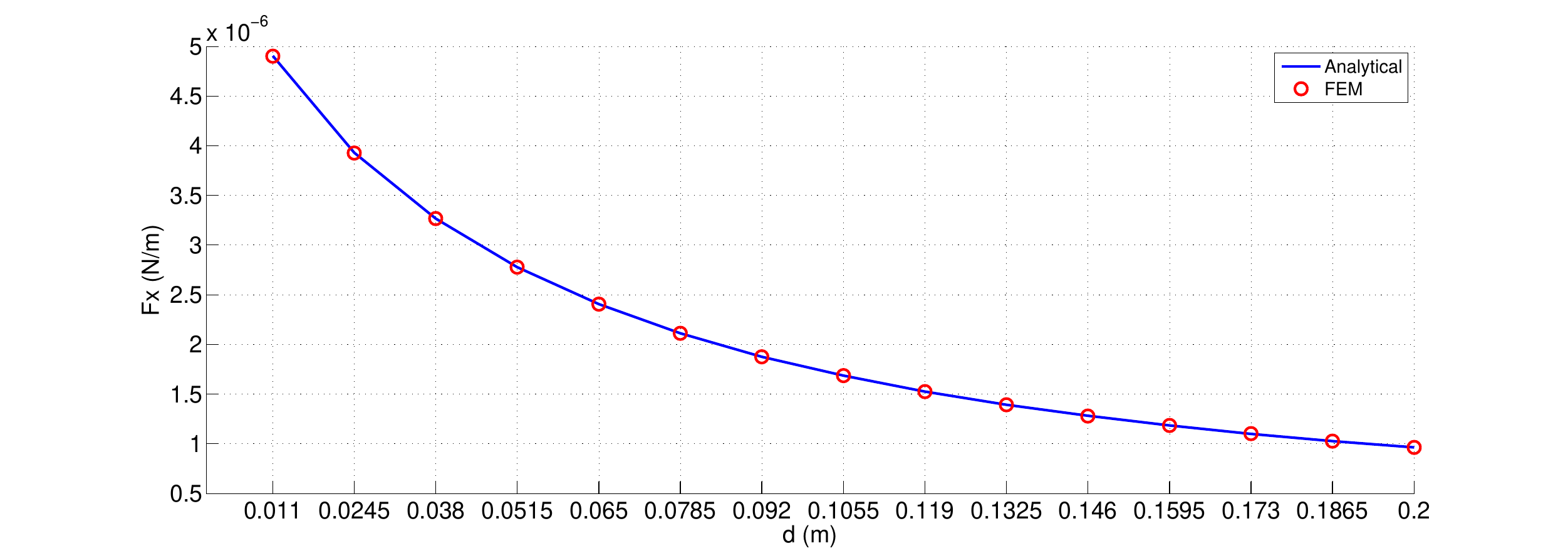}
   \caption{x-component of the force for adjacent massive conductors. $ a~=~0.005\units{m}$, $ b = 0.05\units{m}$}
	  \label{fig:result_example1}
\end{figure}

\noindent
In the first section of the code the geometrical parameters and the distance between the two massive conductors are defined. In the second section the currents are defined. In the third section the computation is done. The results of this example are shown in \figref{fig:result_example1} where the straight line is related to the presented code and the dots are related to the same computation by means of the FEM code.

%

\subsection{Example \#2: adjacent conductors {\small (variable currents)}}
This example describes the use of the function {\tt adjacent.m} when the massive conductors currents are time dependent. The example refers to the following code:

\begin{lstlisting}
%% Section 1: Geometrical parameters
a = 0.005;  % x dimension (m)
b = 0.05;   % y dimension (m)
d = 0.02;   % Distance betwenn adjacent conductors

%%  Section 2: electrical parameters and currents
f = 50;                 % frequency (Hz)
w = 2*pi*f;             % angular frequency (rad/s)
T = 1/f;                % period (s)
t = linspace(0,T,500);  % array with time (500 points)
I = 1;                  % amplitude of the currents
i1 = I*sin(w*t);        % current in conductor 1 (A)
i2 = I*sin(w*t+pi/2);   % current in conductor 2 (A)

%%  Section 3: computation
Fx = adjacent(a,b,d,i1,i2);
\end{lstlisting}

\noindent
In the first section of the code the geometrical parameters and the distance between the two massive conductors are defined. The second section defines the electrical parameters. The frequency is 50 Hz (line 7) and the analysis is performed in time domain considering the discretization of one period (line 9) in 500 points (line 10). This example points out that, for a given geometrical configuration, the function {\tt adjacent.m} is able to handle currents in time domain defined as MATLAB arrays. In other word, it is not necessary to perform a cycle over the time instants.

\begin{figure}[!ht]
   \centering
		\includegraphics[angle=0, width=1\columnwidth]{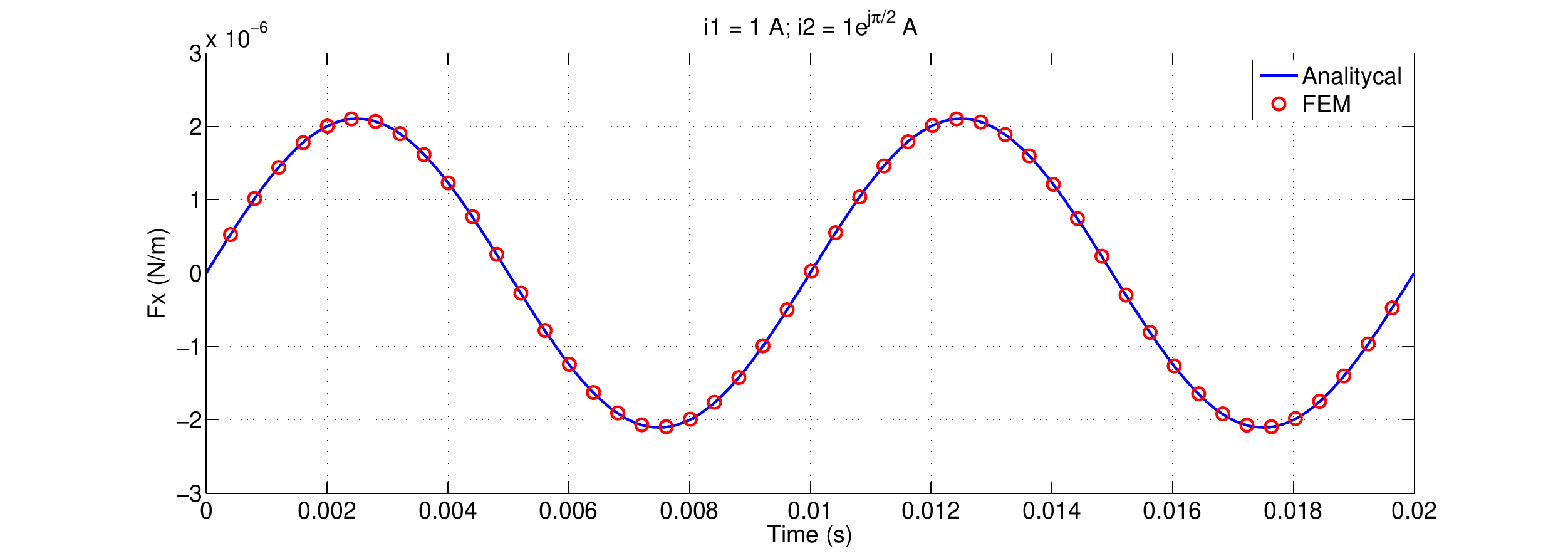}
   \caption{x-component of the force for adjacent massive conductors. $ a~=~0.005\units{m}$, $ b = 0.05\units{m}$}
	  \label{fig:result_example2}
\end{figure}
%
\subsection{Example \#3: non-adjacent conductors}
This example describes the use of the functions {\tt NonAdjacentX.m} and {\tt NonAdjacentY.m}. It is carried out a parametrical analysis by varying the parameters $d$ and $h$.

\noindent
The x-force and y-force can be computed by means of the following code. 

\begin{lstlisting}
%% Section 1: Geometrical parameters
a = 0.005;                  % x dimension (m)
b = 0.05;                   % y dimension (m)
d = linspace(0.011,0.2,15); % x distance between conductors (m)
h = linspace(0.11,0.2,8);   % y distance between conductors (m)

%%  Section 2: conductors currents
i1 = 1; % current in conductor 1 (A)
i2 = 1; % current in conductor 2 (A)

%%  Section 3: computation
for i = 1:length(d)
    for j = 1:length(h)
        Fx(j,i) = NonAdjacentX(a,b,d(i),h(j),i1,i2);
        Fy(j,i) = NonAdjacentY(a,b,d(i),h(j),i1,i2);
    end
end
\end{lstlisting}

\begin{figure}[!ht]
   \centering
		\includegraphics[angle=0, width=1\columnwidth]{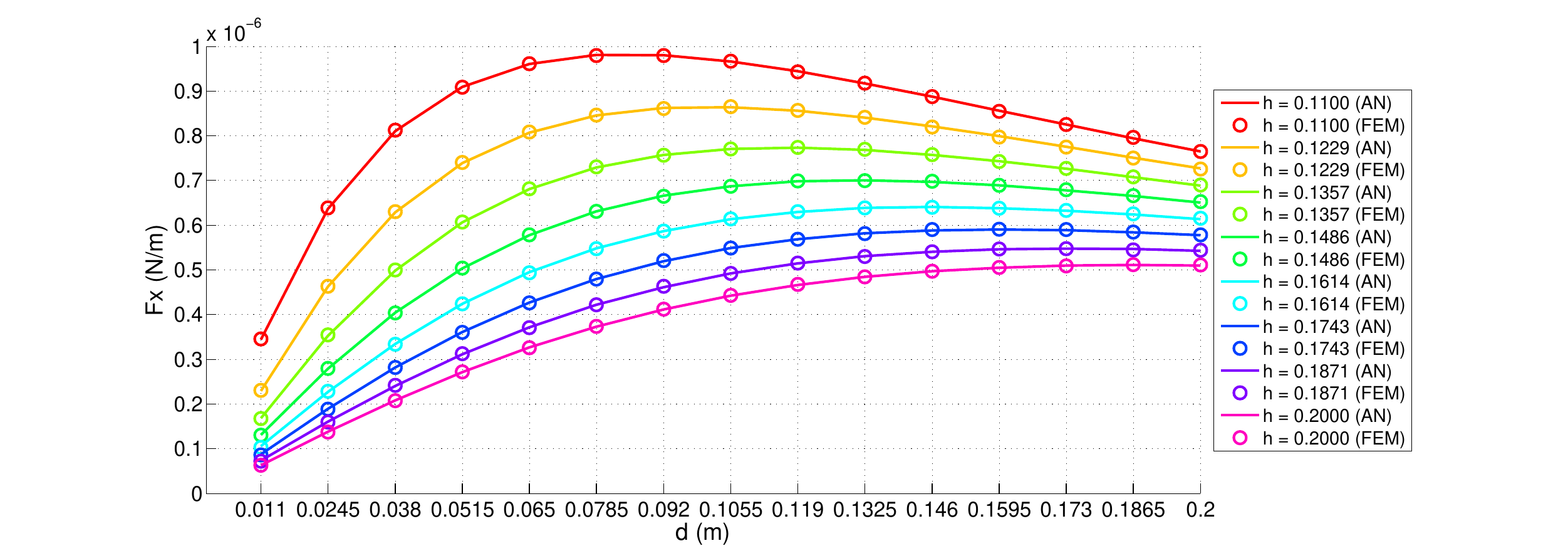}
   \caption{x-component of the force for adjacent massive conductors. $ a~=~0.005\units{m}$, $ b = 0.05\units{m}$}
	  \label{fig:result_example3x}
\end{figure}

\begin{figure}[!ht]
   \centering
		\includegraphics[angle=0, width=1\columnwidth]{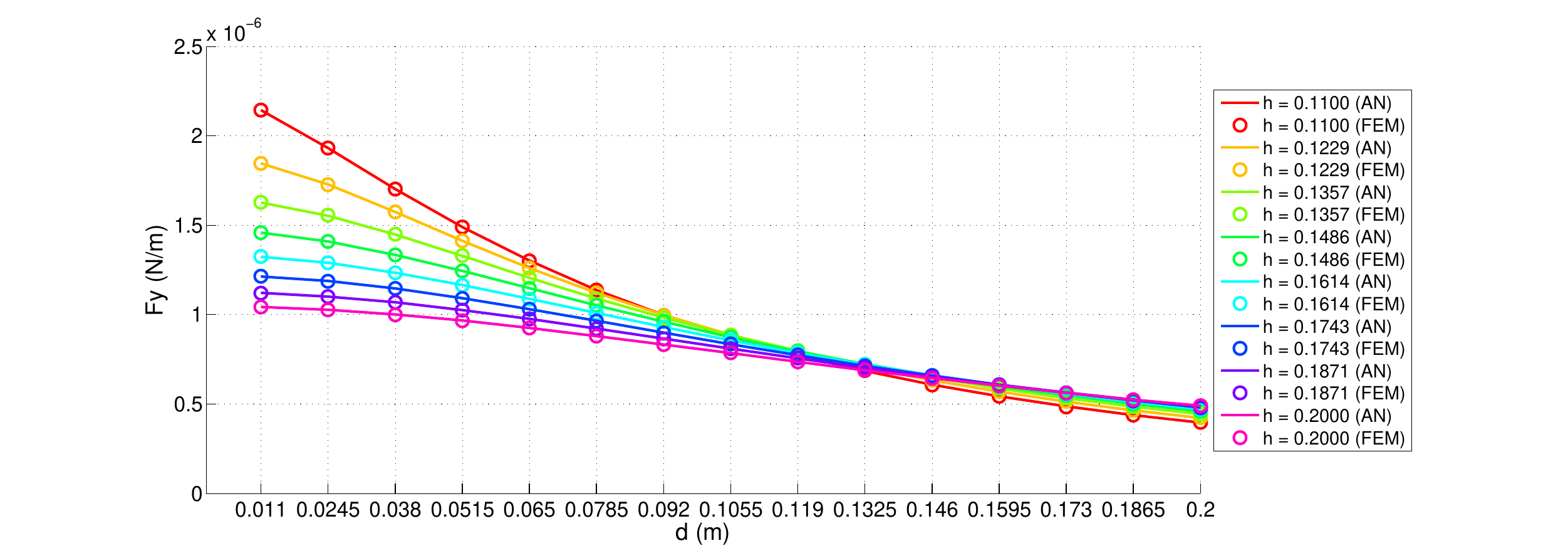}
   \caption{y-component of the force for adjacent massive conductors. $ a~=~0.005\units{m}$, $ b = 0.05\units{m}$}
	  \label{fig:result_example3y}
\end{figure}

\noindent
In the first section of the code the geometrical parameters and the distances between the two massive conductors are defined. In the second section the currents are defined. In the third section the computation is done. The results of this example are shown in \figref{fig:result_example3x} and \figref{fig:result_example3y} where the straight lines are related to the presented code and the dots are related to the same computation by means of the FEM code.


\bibliographystyle{unsrt}
\bibliography{errata}

\end{document}